\documentclass[11pt]{article}
\usepackage{amsfonts,amsmath,amssymb, latexsym,color,epsfig}
\usepackage{tikz,forest}
\usetikzlibrary{shapes}
\date{}

    \newcommand{\F}{\mathcal{F}}

    \newcommand{\HH}{\mathcal{H}}

    \newcommand{\rk}{{\mbox{\rm rk}}}

\def\qed{\ifhmode\unskip\nobreak\fi\quad\ifmmode\Box\else$\Box$\fi}

\title{Exchange properties of finite set-systems}

\author{
{\sl Peter Frankl}\thanks{R\'enyi Institute, P.O.Box 127 Budapest, 1364 Hungary, partially supported by the Ministry of Education and Science of the Russian Federation in the framework of MegaGrant 075-15-2019-1926. Email: \texttt{peter.frankl@gmail.com}.}
\and
{\sl J\'{a}nos Pach}\thanks{R\'enyi Institute, P.O.Box 127 Budapest, 1364 Hungary, partially supported by ERC Advanced Grant GeoScape, NKFI grants KKP-133864 and K-131529, and by the Ministry of Education and Science of the Russian Federation in the framework of MegaGrant 075-15-2019-1926. Email: \texttt{pach@cims.nyu.edu}}
\and
{\sl D\"om\"ot\"or P\'alv\"olgyi}\thanks{MTA-ELTE Lend\"ulet Combinatorial Geometry Research Group, Institute of Mathematics, E\"otv\"os Lor\'and University, Budapest, Hungary. Research supported by the Lend\"ulet program of the Hungarian Academy of Sciences (MTA), under the grant LP2017-19/2017. Email: \texttt{domotorp@gmail.com}.\newline
The extended abstract of this paper has already appeared in the proceedings of the EuroComb 2021 conference.}
}

\date{}

\begin{document}

\maketitle

\begin{abstract}
In a recent breakthrough, Adiprasito, Avvakumov, and Karasev constructed a triangulation of the $n$-dimensional real projective space with a subexponential number of vertices. They reduced the problem to finding a small downward closed set-system $\cal F$ covering an $n$-element ground set which satisfies the following condition:
For any two disjoint members $A, B\in\cal F$, there exist $a\in A$ and $b\in B$ such that
either $B\cup\{a\}\in\cal F$ and $A\cup\{b\}\setminus\{a\}\in\cal F$,
or $A\cup\{b\}\in\cal F$ and $B\cup\{a\}\setminus\{b\}\in\cal F$.
Denoting by $f(n)$ the smallest cardinality of such a family $\cal F$, they proved that $f(n)<2^{O(\sqrt{n}\log n)}$, and they asked for a nontrivial lower bound. It turns out that the construction of Adiprasito {\em et al.}~is not far from optimal; we show that $2^{(1.42+o(1))\sqrt{n}}\le f(n)\le 2^{(1+o(1))\sqrt{2n\log n}}$.
\smallskip

We also study a variant of the above problem, where the condition is strengthened by also requiring that for any two disjoint members $A, B\in\cal F$ with $|A|>|B|$, there exists $a\in A$ such that $B\cup\{a\}\in\cal F$. In this case, we prove that the size of the smallest $\cal F$ satisfying this stronger condition lies between $2^{\Omega(\sqrt{n}\log n)}$ and $2^{O(n\log\log n/\log n)}$.
\end{abstract}

\section{Introduction}\label{intro}

It is an old problem to find a triangulation of the $n$-dimensional real projective space with as few vertices as possible. Recently, Adiprasito, Avvakumov, and Karasev~\cite{AAK} broke the exponential barrier by finding a construction of size $2^{O(\sqrt{n}\log n)}$. For the proof, they considered the following problem in extremal set theory.
\smallskip

What is the {\em minimum} cardinality of a system $\cal F$ of subsets of $[n]=\{1,2,\ldots, n\}$, which satisfies three conditions:
\begin{enumerate}
\item $\cal F$ is {\em atomic}, that is, $\emptyset\in\F$ and $\{a\}\in\cal F$ for every $a\in [n]$; 
\item $\cal F$ is {\em downward closed}, that is, if $A\in\cal F$, then $A'\in\cal F$ for every $A'\subset A$;
\item for any two disjoint members $A, B\in\cal F\setminus \{\emptyset\}$, there exist $a\in A$ and $b\in B$ such that

{\em either} $B\cup\{a\}\in\cal F$ and $A\cup\{b\}\setminus\{a\}\in\cal F$,

{\em or} $A\cup\{b\}\in\cal F$ and $B\cup\{a\}\setminus\{b\}\in\cal F$.
\end{enumerate}
\smallskip

Letting $f(n)$ denote the minimum size of a set-system $\cal F$ with the above three properties, Adiprasito {\em et al.}~proved
\begin{equation}\label{eq1}
f(n)\le 2^{(1/2+o(1))\sqrt{n}\log n},
\end{equation}

where \emph{log} always denotes the base 2 logarithm.
They used the following construction. Let $s,t>0$ be integers, $n=st$. Fix a partition $[n]=X_1\cup\ldots\cup X_t$ of the ground set into $t$ parts of equal size, $|X_1|=\ldots=|X_t|=s$. Let
\begin{equation}\label{eq1.5}
\F=\cup_{i=1}^t{\F}_i,\;\;\mbox{ where }\;\;
{\cal F}_i=\{ F\subseteq [n] : |F\cap X_j|\le 1 \mbox{ for every } j\not= i\},
\end{equation}

for $1\le i\le t$. (In the definition of ${\cal F}_i$, there is no restriction on the size of $F\cap X_i$.)
It is easy to verify that ${\cal F}$ meets the requirements. We have
$$|{\cal F}|=(t2^s-(s+1)(t-1))(s+1)^{t-1}<2^{s+t\log(s+1)+\log t}.$$

Substituting $s=(1/\sqrt2+o(1))\sqrt{n\log n}$ and $t=(\sqrt2+o(1))\sqrt{n/\log n}$, we obtain that
\begin{equation}\label{eq2}
f(n)\le 2^{(1/\sqrt2+o(1))\sqrt{n\log n}+(\sqrt2+o(1))\sqrt{n/\log n}\cdot\log\sqrt n}
=2^{(1+o(1))\sqrt{2n\log n}}.
\end{equation}

This is slightly better than (\ref{eq1}). (The authors of \cite{AAK} remarked that their bound can be improved by a ``subpolynomial factor.'') Any further improvement on the upper bound would result in a smaller triangulation of the $(n-1)$-dimensional projective space.
\smallskip

Our first theorem implies that (\ref{eq2}) is not far from optimal.

The \emph{rank} of a set-system $\cal F$, denoted by $\rk({\cal F})$, is the size of the largest set $F\in\F$; see, {\em e.g.}, \cite{D}.

We denote by $\lfloor x\rceil$ the integer closest to $x$. (We will use this notation only for $x=\sqrt{2n}$, in which case $\lfloor x\rceil$ is uniquely determined.) 
\medskip

\noindent{\bf Theorem 1.} {\em Let $\cal F$ be an atomic system of subsets of $[n]$, such that for any two disjoint members $A, B\in\cal F$, {\em either} there exists $a\in A$ such that $B\cup\{a\}\in\cal F$, {\em or} there exists $b\in B$ such that $A\cup\{b\}\in\cal F$. Then we have

{\bf (i)}\;\;\;  $|\F|\ge e^{(2e^{-1/\sqrt 2}+o(1))\sqrt{n}}\ge 2^{(1.42+o(1))\sqrt{n}}$; 

{\bf (ii)}\;\; $\rk(\F)\ge \lfloor\sqrt{2n}\rceil$, and this bound is best possible.}
\medskip

\noindent
If we also assume that $\F$ is downward closed, then the inequality $\rk(\F)\ge \lfloor\sqrt{2n}\rceil$ immediately implies that $|\F|\ge 2^{\lfloor\sqrt{2n}\rceil}$. This is only slightly weaker than the lower bound $f(n)\ge 2^{(1.42+o(1))\sqrt{n}}$, which follows from part (i).

\medskip

\noindent{\bf Remark.}
We remark that the assumptions of Theorem 1 are weaker than those made by Adiprasito {\em et al.}, in two different ways: we do not require that $\cal F$ is downward closed (which is their condition 2), and the exchange condition between two disjoint sets is also less restrictive than condition 3.
Nevertheless, we know no significantly smaller set-systems satisfying these weaker conditions than the ones described in (\ref{eq1.5}), for which $|{\cal F}|=2^{(1+o(1))\sqrt{2n\log n}}$.
We can, however, further weaken the conditions under which Theorem 1 holds; instead of the exchange property, it is sufficient to assume the following:\\
\indent
{\em For any two disjoint members $A,B\in\F$ with $|A|=|B|$, there is a set $C\in\F$ such that $C\subset A\cup B$ and $|C|=|A|+1$.}\\
\indent
This answers Question 3.7 of Adiprasito {\em et al.}~\cite{AAK}:
from Claim 3.1 of \cite{AAK}, one cannot obtain a significantly better construction, using another family. To see that condition (3) in Claim 3.1 implies our condition above, apply it with the unit vector $X$ identified with $A\cup B$.

\smallskip

\medskip

While part (ii) of Theorem 1 is tight, we suspect that part (i) and the lower bound $f(n)\ge 2^{\Omega(\sqrt{n})}$ can be improved. As a first step, we slightly strengthen the assumptions of Theorem 1, in order to obtain a better lower bound on $|\cal F|$.

\medskip

\noindent{\bf Theorem 2.} {\em Let $\cal F$ be an atomic system of subsets of $[n]$, such that for any two disjoint members $A, B\in\cal F$, {\em either} there exists $a\in A$ such that $B\cup\{a\}\in\cal F$, {\em or} there exists $b\in B$ such that $A\cup\{b\}\in\cal F$. Moreover, suppose that if $|A|<|B|$, then the second option is true.

Then we have $|\F|\ge 2^{(1/2+o(1))\sqrt{n}\log n}$.}
\medskip

This lower bound exceeds the upper bound in (\ref{eq2}). Therefore, construction (\ref{eq1.5}) cannot satisfy the stronger assumptions in Theorem 2.
For example, set
$$A=\{a_1\}\cup\{a_2,a_2'\}\cup\{\emptyset\}\cup\ldots\cup\{\emptyset\}\in{\cal F}_2\subset\F,$$ $$B=(X_1\setminus\{a_1\})\cup\emptyset\cup\{\emptyset\}\cup\ldots\cup\{\emptyset\}\in{\cal F}_1\subset \F,$$
where $a_1\in X_1$ and $a_2,a_2'\in X_2$. If $s>4$, then $|A|<|B|$, but there is no element of $B$ that can be added to $A$ such that the resulting set also belongs to $\cal F$.
If $s\le 4$, then the conditions of Theorem 2 are satisfied, but the construction is uninteresting, as $|\F|=2^{\Theta(n)}$ and $\rk(\F)=\Theta(n)$.
\smallskip

Our next result provides a nontrivial construction.

\medskip
\noindent{\bf Theorem 3.} {\em There exists an atomic downward closed set-system ${\cal F} \subset 2^{[n]}$ with the property that for any two disjoint members $A, B\in\cal F$ with $|A|\le|B|$, there is $b\in B$ such that $A\cup\{b\}\in\cal F$, and

{\bf (i)}\;\;\; $|\F|\le 2^{(2+o(1))n\log\log n/\log n},$

{\bf (ii)}\;\;  $\rk(\F)\le (2+o(1))n/\log n$.}
\medskip

The proofs of Theorems 1, 2, and 3 are presented in Sections 2, 3, and 4, respectively.

\section{Proof of Theorem 1}\label{sec2}

We start with a statement which immediately implies the inequality in part (ii).
\medskip

\noindent{\bf Lemma 2.1.} {\em Let $k\ge 1$ be an integer, $n>{k\choose 2}$, and let $\F$ be an atomic family of subsets of $[n]$ satisfying the exchange property in Theorem 1, or the condition in the Remark.

Then there is a set $F\in\F$ such that $|F|=k$. This bound cannot be improved: there are families satisfying the conditions, for which $\rk({\F})= k$. }
\medskip

\noindent{\bf Proof.} By induction on $k$. For $k=1$, the claim is trivial.
Suppose that $k>1$ and that the lemma has already been proved for $k-1$.

Let $\F\subset 2^{[n]}$ be a family satisfying the conditions, where $n > \binom{k+1}{2}$. By the induction hypothesis, there is a member $A\in\F$ such that $|A|=k$. 
Consider the family $\F'=\{F\in\F : F\cap A=\emptyset\}$. Obviously, $\F'$ satisfies the conditions on the ground set $[n]\setminus A$, and we have
$|[n]\setminus A|>\binom{k+1}{2}-k=\binom{k}{2}$. Hence, we can apply the induction hypothesis to $\F'$ to find a set $B\in \F'$ of size $k$ which is disjoint from $A$. Using the exchange property in Theorem 1, or the condition in the Remark, for the sets $A$ and $B$, we can conclude that $\F$ has a member with $k+1$ elements. 

Now we show the tightness of Lemma 2.1.
Let $X_1, \ldots, X_{k-1}$ be pairwise disjoint sets with $|X_i|=i$, for every $i$. Then $V=X_1\cup\ldots\cup X_{k-1}$ is a set of ${k\choose 2}$ elements. For $i=1,\ldots,k-1$, define
\begin{equation}\label{eq2.5}
\F_i=\{ F\subseteq V : |F\cap X_j|=0 \mbox{ for every } j<i \mbox{ and } |F\cap X_j|\le 1 \mbox{ for every } j>i\}.
\end{equation}
In the definition of $\F_i$, there is no restriction on the size of $F\cap X_i$. Let $\F=\F_1\cup\ldots\cup\F_{k-1}.$  Obviously, every member of $\F_i$ has at most $|X_i|+k-1-i=k-1$ elements, which yields that $\rk(\F)=\max_{i=1}^k\rk(\F_i)= k-1$. Furthermore, $\F$ is atomic and any two disjoint members of $\F$ satisfy the exchange condition in Theorem 1. Hence, the lemma is tight. \hfill $\Box$
\medskip

We remark that the maximal sets in the above $\F$ form the same hypergraph as the one defined in Example 3 of \cite{L} for $v=1$. 
\smallskip

To prove the inequality $\rk(\F)\ge\lfloor\sqrt{2n}\rceil$ in part (ii) of Theorem 1, we have to find the largest $k$ for which we can apply Lemma 2.1. It is easy to verify by direct computation that
$$\max\{k : \binom{k}{2}<n\}=\lfloor\sqrt{2n}\rceil.$$

If $n={k\choose 2}$ for some $k\ge 1$, then the tightness of part (ii) of Theorem 1 follows from the tightness of Lemma 2.1. Suppose next that ${k\choose 2}<n<{k+1\choose 2}$.
Let $X_1,\ldots, X_k$ be pairwise disjoint sets with $|X_i|=i$ for every $i<k$ and let $|X_k|=n-{k\choose 2}$. Set $V=X_1\cup\ldots\cup X_k.$ For $i=1,\ldots, k,$ define $\F_i$ as in (\ref{eq2.5}), and let
$\F=\F_1\cup\ldots\cup\F_k$. Then $\F$ has the exchange property and $\rk(\F)=k=\lfloor \sqrt{2n}\rceil.$ This proves part (ii) of Theorem 1.
\medskip

It remains to establish part (i). Let $\F$ be a family of subsets of $[n]$ satisfying the conditions. 
Let $\F'$ denote the {\em $k$-uniform hypergraph} (i.e., family of $k$-element sets) consisting of all $k$-element sets in $\F$, i.e., $\F'=\{F\in\F : |F|= k\}$.

The {\em independence number} $\alpha(\HH)$ of a hypergraph $\HH$ is the maximum cardinality of a subset of its ground set which contains no element (hyperedge) of $\HH$.
It follows from Lemma 2.1 that any subset $S\subseteq [n]$ of size $|S|=\binom k2+1$ contains at least one element of $\F$ whose size is $k$. Therefore, any such set contains at least one element of $\F'$, which means that $\alpha(\F')\le \binom k2$.

We need a result of Katona, Nemetz, and Simonovits \cite{KNS} which is a generalization of Tur\'an's theorem to $k$-uniform hypergraphs.
\medskip

\noindent{\bf Lemma 2.2.} \cite{KNS} {\em Let $\HH$ be a $k$-uniform hypergraph on an $n$-element ground set. If the independence number of $\HH$ is at most $\alpha$, then we have $$|\HH|\ge\binom nk\bigg/\binom \alpha k.$$}
\medskip

Applying Lemma 2.2 to the hypergraph $\HH=\F'$ with $k=(\sqrt{2}e^{-1/\sqrt 2}+o(1))\sqrt{n}$ and  $\alpha=\binom k2$, we obtain
$$|\F|\ge |\F'|\ge e^{(2e^{-1/\sqrt 2}+o(1))\sqrt{n}}\ge 2^{(1.42+o(1))\sqrt n},$$
completing the proof of part (i). This bound is slightly better than the inequality $|\F|\ge 2^{\lfloor\sqrt{2n}\rceil}$, which immediately follows from part (ii), under the stronger assumption that $\F$ is downward closed.

\section{Proof of Theorem 2}\label{sec3}

Let $\cal F$ be an atomic set-system on an $n$-element ground set $X$, where $n$ is large, and let $s$ and $t$ be two positive integers to be specified later.
We describe a procedure to identify $\sum_{i=0}^ts^i$ distinct members of $\cal F$.
To explain this procedure, we fix an $s$-ary tree $T$ of depth $t$.
At the end, each of the $s^t$ root-to-leaf paths in $T$ will correspond to a unique member of $\F$.
\smallskip

Each {\em non-leaf vertex} $v$ will be associated with an $s$-element subset $X(v)\subset X$ such that along every root-to-leaf path $p=v_0v_1\ldots v_t$, the sets $X(v_0),X(v_1),\ldots,X(v_{t-1})$, associated with the root and with the internal vertices of $p$, will be pairwise disjoint. See Figure 1 for an example.

Each {\em edge} $e=vu$ of $T$, where $u$ is a child of $v$, will be labelled with an element $x(e)\in X(v)$ in such a way that every edge from $v$ to one of its $s$ children gets a different label. Thus,
$$\{x(vu) : u \mbox{ is a child of }v\}=X(v).$$

Denoting the root by $v_0$, we choose $X(v_0)$ to be an arbitrary $s$-element subset of the ground set $X$, and set $F(v_0)=\emptyset\in\F$. For any non-root vertex $v$, let $$F(v)=\{x(e) : e \mbox{ is an edge along the root-to-}v\mbox{ path}\}.$$

We will choose $X(v)$ such that $F(v)\in\F$ for every $v$. All of the sets $F(v)$ will be distinct, as any two different paths starting from the root diverge somewhere, unless one contains the other.

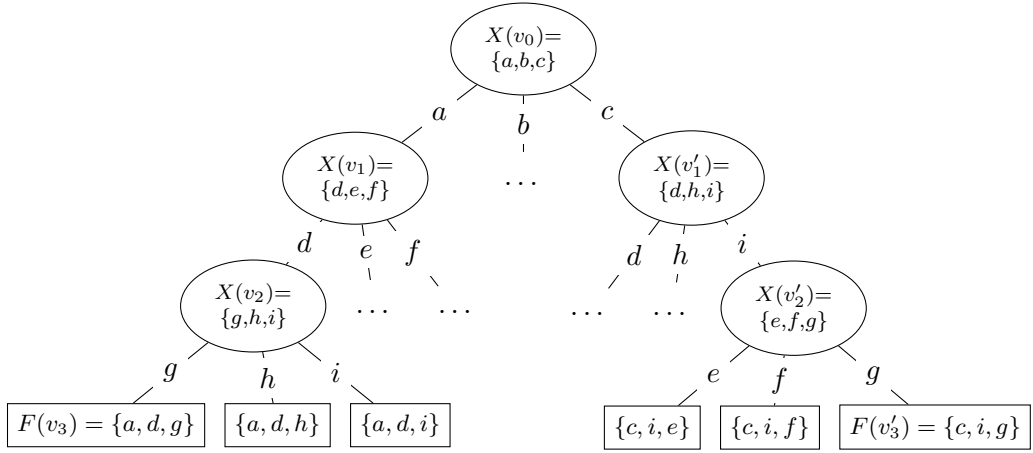
\begin{figure}[h]
	\begin{center}
		\begin{forest}
			csucs/.style={ellipse, draw},
			level/.style={rectangle,draw, before packing={l+=5mm},font=\footnotesize},
			pontok/.style={circle,draw=white}
			[{$X(v_0)=\atop \{a,b,c\}$}, for tree=csucs
			[{$X(v_1)=\atop  \{d,e,f\}$}, edge label={node[midway,fill=white]{$a$}}
			[{$X(v_2)=\atop  \{g,h,i\}$}, edge label={node[midway,fill=white]{$d$}}
			[{$F(v_3)=\{a,d,g\}$}, level, edge label={node[midway,fill=white]{$g$}}
			]
			[{$\{a,d,h\}$}, level, edge label={node[midway,fill=white]{$h$}}
			]
			[{$\{a,d,i\}$}, level, edge label={node[midway,fill=white]{$i$}}
			]
			]		
			[{$\dots$}, pontok, edge label={node[midway,fill=white]{$e$}}
			]
			[{$\dots$}, pontok, edge label={node[midway,fill=white]{$f$}}
			]
			]
			[,phantom]
			[{$\dots$}, pontok, edge label={node[midway,fill=white]{$b$}}]
			[,phantom]
			[{$X(v_1')=\atop  \{d,h,i\}$}, edge label={node[midway,fill=white]{$c$}}		
			[{$\dots$}, pontok, edge label={node[midway,fill=white]{$d$}}
			]
			[{$\dots$}, pontok, edge label={node[midway,fill=white]{$h$}}
			]		
			[{$X(v_2')=\atop  \{e,f,g\}$}, edge label={node[midway,fill=white]{$i$}}
			[{$\{c,i,e\}$}, level, edge label={node[midway,fill=white]{$e$}}
			]
			[{$\{c,i,f\}$}, level, edge label={node[midway,fill=white]{$f$}}
			]		
			[{$F(v_3')=\{c,i,g\}$}, level, edge label={node[midway,fill=white]{$g$}}
			]	
			]
			]	
			]
		\end{forest}
		\caption{Construction of the auxiliary tree $T$ for the proof of Theorem 2.}
	\end{center}
\end{figure}

Suppose that we have already determined the set $X(u)$ for all ancestors of some non-leaf vertex $v$ at level $\ell<t$ of $T$. At this point, we already know the set $F(v)\in\F$, where $|F(v)|=\ell$, and we want to determine $X(v)$. The next lemma guarantees that there is a good choice for $X(v)$.

\medskip

\noindent{\bf Lemma 3.1.} {\em There is an $s$-element subset $X(v)\subset X$ such that for every $x\in X(v)$, we have $F(v)\cup \{x\}\in\F$.}
\smallskip

\noindent{\bf Proof.} Let $Z=\cup \{X(u) : u\mbox{ lies on the root-to-}v\mbox{ path}, u\neq v\}$. If $v$ is at level $\ell<t$, we have $|Z|=\ell s\le (t-1)s$.
Let $$Y=\{y\in X\setminus Z : F(v)\cup \{y\}\notin \F\}.$$
In other words, $Y$ consists of all elements of $X\setminus Z$ that cannot be added to $F(v)$ to obtain a set in $\F$.

Consider the family $\F'=\{F\in\F : F\subseteq Y\}$.
If $|Y|> \binom{t}{2}$, then Lemma 2.1 implies that 
there is a set $B\in \F'$ with $|B|=t>|F(v)|$. In this case, we can apply the exchange condition in Theorem 2 to the sets $F(v)$ and $B$, to conclude that there exists $b\in B$ for which $F(v)\cup\{b\}\in \F$. However, this contradicts the fact that $b\in Y$.

Thus, we can assume that $|Y|\le \binom{t}{2}$. Now we have
$$|(X\setminus Z)\setminus Y|\ge n-(t-1)s-{t\choose 2}.$$
If the right-hand side of this inequality is at least $s$, there is a proper choice for the set $X(v)$.
For this, it is enough if $n \ge \frac{t^2}2+ts$, or, equivalently, $2 \ge (\frac{t}{\sqrt n})^2+\frac{t}{\sqrt n}\frac{2s}{\sqrt n}$.
To achieve this, let $n$ be large,
$s=\lfloor \sqrt{n}/\log^2 n\rfloor$, and $t=\lfloor(1-1/\log n)\sqrt{2n}\rfloor$.
$\Box$
\medskip

By the above procedure, we can recursively assign a different set $F(v)\in\F$ to each vertex $v$ of $T$. This gives
\[|\F|\ge \sum_{i=0}^t s^i\ge (n^{1/2}/\log^2 n)^{(1-o(1))\sqrt{2n}}=n^{\sqrt{(1/2+o(1))n}},
\]
as desired.

\section{Proof of Theorem 3}
Assume for simplicity that $n$ is a multiple of $k$, and fix a partition $[n]=X_1\cup\ldots\cup X_{n/k}$ into $n/k$ parts, each of size $k$. That is, let $|X_1|=\dots=|X_{n/k}|=k$, where $k$ is the largest number for which $2^{k-2}\le n/k$; this gives $k=(1+o(1))\log n$.
We will also assume $n,k\ge 3$.
\smallskip

For any $A\subset [n]$ and $0\le i\le k$, let $p_A(i)$ and $s_A(i)$ denote the number of parts $X_t$ which intersect $A$ in {\em precisely} $i$ elements and in {\em at least} $i$ elements, respectively. Thus, we have $s_A(i)=\sum_{j=i}^k p_A(j)$ and $|A|=\sum_{i=1}^k ip_A(i)=\sum_{i=1}^ks_A(i).$  Define the {\em profile vector} of $A$, as $$p_A=(p_A(k), p_A(k-1), \ldots,p_A(0)),$$ and let $$s_A=(s_A(k), s_A(k-1), \ldots,s_A(0)).$$

That is, $p_A(0)$ is the number of parts that are disjoint from $A$, while $s_A(0)$ is always equal to $n/k$. We claim that the set-system
$${\cal F}=\{ A\subseteq [n]: s_A(k)\le 1 \mbox{ and } s_A(i)\le 2^{k-1-i} \mbox{ for every } 2\le i\le k-1 \} $$
meets the requirements of the theorem.
Notice that for $i=0$ and $i=1$, there is no restriction on $s_A(i)$ other than the trivial bounds $0\le s_A(i)\le n/k$.
In particular, if $A$ is an element of $\cal F$ with maximum cardinality, we have
$$s_A=(1,1,2^1,2^2,2^3,\ldots,2^{k-3},n/k,n/k),\; \mbox{ and}$$ $$p_A=(1,0,1,2^1,2^2,\ldots,2^{k-4},n/k-2^{k-3},0).$$
Here, we used that $2^{k-3}\le n/k$, by assumption.
Thus, the size of such a largest set is
\[
\rk(\F)=|A|=\sum_{i=1}^k s_A(i)=n/k+\sum_{i=1}^{k-1} 2^{k-i-1}+1=n/k+2^{k-2}\le 2n/k=(2+o(1))n/\log n.
\]
This also gives the following simple bound on the size of the family.
\[
|\F|\le \sum_{i=0}^{\rk(\F)} \binom{n}{i}
\le \sum_{i=0}^{2n/k} \binom{n}{i}\le \left({en\over 2n/k}\right)^{2n/k}
\!\le \left({e\log n\over 2}\right)^{(2+o(1))n/\log n}
\!\!\le 2^{(2+o(1))n\log\log n/\log n}.
\]
\smallskip

To prove that $\cal F$ meets the requirements of the theorem, we consider any pair of disjoint sets $A,B\in \cal F$ such that there exists no $b\in B$ such that $A\cup\{b\}\in\cal F$. We need to show that in this case we have $|A|>|B|$.
Suppose for contradiction that this is not the case, and fix a {\em counterexample} for which $|B|\ge|A|$ and $|B|-|A|$ is as large as possible. We refer to such a counterexample as a {\em maximal counterexample.}
We say that $s_A(i)$ is \emph{saturated} if $s_A(i)=2^{k-1-i}$ for $1<i<k$, if $s_A(i)=1$ for $i=k$, and if $s_A(i)=n/k$ for $i=1$.
\medskip

\noindent{\bf Claim 4.1.}  {\em For any part $X\in\{X_1,\ldots,X_{n/k}\},$ if $B\cap X\not=\emptyset$, then $A\cap X\not=\emptyset$.
Moreover, if $|A\cap X|=i$, then $s_A(i+1)$ is saturated.
Hence, $s_B(k)=0$.}
\medskip

\noindent{\bf Proof.} Otherwise, we could add one element from $B\cap X$ to $A$, as $s_A(i+1)$ was not saturated. As $A$ and $B$ are disjoint, this also implies $s_B(k)=0$. \hfill  $\Box$
\medskip

\noindent{\bf Claim 4.2.}  {\em There is a maximal counterexample $A, B\in\cal F$ such that, for every $X\in\{X_1,\ldots,X_{n/k}\},$

{\bf (i)}\;\;\; $A\cap X\not=\emptyset$;

{\bf (ii)}\;\; if $|A\cap X|=1$, then $B\cap X\not=\emptyset.$}
\medskip

\noindent{\bf Proof.}
 Among all maximal counterexamples $A,B\in\cal F$, choose one for which $s_A(1)$ is as large as possible.
Let $j\ge 1$ be the smallest positive integer for which there is a part $X$ with $|A\cap X|=j$.
Thus, we have $s_A(j)=s_A(j-1)=\ldots=s_A(1)\ge s_B(1)$, where the last inequality follows from Claim 4.1.
This implies that $s_A(i)$ cannot be saturated for any positive integer $i<j$.

Suppose first that $j>k/2$. Then we have
$$|A|=\sum_{i=1}^{n/k} |A\cap X_i|\;
>\sum_{X_i\cap A\ne \emptyset} k/2\;
\ge \sum_{\substack{X_i\cap A\ne \emptyset\\X_i\cap B\ne \emptyset}} k/2\;
\ge \sum_{\substack{X_i\cap A\ne \emptyset\\X_i\cap B\ne \emptyset}} |B\cap X_i|\;
=\sum_{i=1}^{n/k} |B\cap X_i|
=|B|,$$
where the second inequality follows from Claim 4.1, the third inequality from the disjointness of $A$ and $B$ and $j>k/2=|X_i|/2$, and the final equality again from Claim 4.1.
Hence, in this case, the pair $A,B$ did not constitute a counterexample.

From now on, we can assume $j\le k/2$.
If for any part $X$, we have $|A\cap X|=j$ and $B\cap X=\emptyset$, then we claim that $A'=A\setminus X$ and $B$ would also form a counterexample.
To see this, assume for a contradiction that some $b$ can be added to $A'$.
Then $b$ could also be added to $A$, unless $b\in X'$ for some part $X'$ such that $|A'\cap X'|=j-1$.
But this is impossible for $j\ge 2$, because of the minimality of $j$. It is also impossible for $j=1$, because, by Claim 4.1, we have $A'\cap X'\ne j-1$.
For the counterexample formed by $A'$ and $B$, we have $|B|-|A'|>|B|-|A|$, contradicting the maximality of $|B|-|A|$.
This proves (ii) because $|A\cap X|=1$ implies $j=1$ by the definition of $j$.
To prove (i), note that
by the `moreover' part of Claim 4.1, we can conclude that $s_A(j+1)$ is saturated.

Next, we show that $s_A(j)$ is also saturated. 
Otherwise, if (i) does not hold, pick a part $X$ for which $A\cap X=\emptyset$.
By Claim 4.1, this implies that $B\cap X=\emptyset$.
Add any $j\le k/2$ elements of $X$ to $A$ and $j\le k/2$ other elements of $X$ to $B$, and denote the resulting sets by $A'$ and $B'$, so that $|B'|-|A'|=|B|-|A|$. In view of Claim 4.1, we have $s_{B'}(j) = s_{B}(j)+1\le s_{B}(1)+1\le s_{A}(1)+1=s_A(j)+1=s_{A'}(j).$ Thus, $A'$ and $B'$ belong to $\cal F$, as we assumed that $s_A(j)$ is not saturated. To show that they constitute a counterexample, we need to prove that for every $b\in B'$, we have $A'\cup\{b\}\notin\F$.
That is, if $b\in X'$ for some part $X'$ and $|A'\cap X'|=i$, then we need to prove that $s_{A'}(i+1)$ is saturated.
We know that $i\ge j$ as the intersection of $A'$ with any part is at least $j$.
If $i>j$, then $X'\ne X$ and so we could also add $b$ to $A$ from $B$ to obtain $A\cup\{b\}\in\F$, contradicting that they formed a counterexample.
If $i=j$, then $s_{A'}(j+1)=s_A(j+1)$ is saturated because of the previous paragraph.
Since $s_{A'}(1)>s_A(1)$, this contradicts the maximal choice of $A$ made at the very beginning of this proof.

Moreover, $j=1$ must hold. Otherwise, if (i) does not hold, pick a part $X$ disjoint from $A$ and $B$, and add any $j-1$ elements of $X$ to $A$ and $j-1$ other elements of $X$ to $B$. Denote the resulting sets by $A'$ and $B'$, so that $|B'|-|A'|=|B|-|A|$. Similarly as before, these new sets also belong to $\F$, as $s_A(j)$ was saturated.
We get a contradiction again, because we have $s_{A'}(1)>s_A(1)$.
Since $s_A(1)$ is saturated, we have $s_A(1)=n/k$. This proves part (i).\hfill $\Box$
\medskip


\noindent{\bf Claim 4.3.} {\em There is a maximal counterexample $A,B\in\cal F$ for which Claim $4.2$ holds and $|B\cap X|>1$ implies $|A\cap X|=1$, for every $X\in\{X_1,\ldots,X_{n/k}\}.$}
\medskip

\noindent{\bf Proof.}
There is a maximal counterexample $A,B\in\cal F$ for which Claim $4.2$ holds.
By definition, we have
$$s_B(2)\le 2^{k-3}\le n/k-2^{k-3}\le n/k-s_A(2)=p_A(0)+p_A(1)=p_A(1),$$
where the second inequality follows from our assumption $n/k\ge 2^{k-2}$.

Therefore, if $|B\cap X|>1$ and $|A\cap X|>1$ for some part $X$, then there exists another part $X'$ for which $|B\cap X'|=1$ and $|A\cap X'|\le 1$.
By Claim 4.2 (i) (or by Claim 4.1), $|A\cap X'|=1$. Choose any $|B\cap X|-1$ elements of $B\cap X$, and remove them from $B$. Choose the same number of elements of $X'\setminus A$, and add them to $B$. By a repeated application of this procedure, we can achieve that the condition in the claim is satisfied. $\Box$
\smallskip

From now on, we consider a counterexample $A,B\in \cal F$ satisfying the condition in Claim 4.3. Let $j$ be the smallest positive integer with the property that for every part $X$ with $|A\cap X|=j$, we have $B\cap X=\emptyset$.
To see that there exists at least one such $j$, notice that if for some $j$ there is no part $X$ with $|A\cap X|=j$, then $j$ meets the requirement. This implies that the number $k-1$ has the desired property: If there is a part $X$ such that $|A\cap X|=k-1$, then since $s_A(k-1)=1$, there can be no part $X'$ such that $|A\cap X'|=k$, so if $B\cap X\neq\emptyset$, we could add $B\cap X$ to $A$, contradicting that they are a counterexample. It follows from Claim 4.2 (ii) that $j=1$ is not possible, so $1<j<k$.

For each part $X$ intersecting $A$ in more than $j$ elements, remove $|A\cap X|-j$ elements of $A\cap X$ from $A$ and remove all elements in $B\cap X$ from $B$.
Denote the sets obtained this way by $A'$ and $B'$, respectively.
As $s_{A'}(i)=s_A(i)$ for every $i\le j$, and if $A'\cap X=j$, then $B\cap X=\emptyset$, the sets $A'$ and $B'$ are a counterexample.
According to Claim 4.3, $|B'|-|A'|\ge |B|-|A|$, and it is easy to see that Claims 4.2 and 4.3 remain valid.
We summarize these properties in the below claim.
\medskip

\noindent{\bf Claim 4.4.} {\em There exists a maximal counterexample $A,B$ which satisfies Claims 4.2 and 4.3, and the following three properties:

{\bf (i)}\;\;\;\;  $s_A(j+1)=0$;

{\bf (ii)}\;\;\;  $p_B(0)\ge s_A(j)$;

{\bf (iii)}\;\;  $s_A(i)=2^{k-1-i}$ for every $2\le i\le j$.}
\medskip

\noindent{\bf Proof.} Parts (i) and (ii) follow directly from the construction, while (iii) follows from the `moreover' part of Claim 4.1 using $j<k$.\hfill  $\Box$
\medskip

Now we can easily complete the proof of Theorem 3. We have

$$|A|=\sum_{i=1}^k s_A(i)=\sum_{i=2}^{j} 2^{k-1-i} + n/k.$$
On the other hand, $s_B(k)=0$ holds, by Claim 4.1, and
$s_B(1)\le n/k-s_A(j)=n/k-2^{k-1-j}$, by Claim 4.4 (iii). Thus,
$$|B|=\sum_{i=1}^k s_B(i)\le n/k-2^{k-1-j}+ \sum_{i=2}^{k-1} 2^{k-1-i},$$
$$|A|-|B|\ge 2^{k-1-j}-\sum_{i=j+1}^{k-1} 2^{k-1-i}=1.$$
This means that $|A|>|B|$, contradicting our assumption that $A,B\in\cal F$ is a counterexample.

\section{Concluding remarks}

If we strengthen the condition of our results by requiring that for any two non-empty disjoint members $A, B\in\cal F$, there exist $a\in A$ and $b\in B$ such that $B\cup\{a\}\in\cal F$ and $A\cup\{b\}\in\cal F$ {\em both} hold, then the problem becomes trivial. Any atomic set-system $\F\subset 2^{[n]}$ with this property must contain all subsets of $[n]$. Indeed, every set $F=\{x_1,\ldots,x_k\}$ can be built up, sequentially applying the condition to the sets $\{x_1,\ldots,x_i\}$ and $\{x_{i+1}\}$, for $i=1,\ldots,k-1$.
\smallskip

In Theorems 1 and 2, we only assume that $\F$ is atomic. However, our best constructions have the stronger property that $\F$ is downward closed. Could we substantially strengthen these results under the stronger assumption? The proof of the bound $|{\F}|\ge 2^{\lfloor\sqrt{2n}\rceil}$, which is only slightly weaker than Theorem 1 (i), becomes much easier if we assume that $\F$ is downward closed, and the proof of Theorem 2 can also be simplified if $\F$ is downward closed.
\smallskip

The property of the set-system described in Theorems 2 and 3 is reminiscent of the \emph{independent set exchange property} of matroids; see \cite{Ox}.
A common generalization of these two properties would be to require that for any two members $A, B\in\cal F$, if either $|A|=|B|$ and $A\cap B=\emptyset$, or $|A|<|B|$ (but they are not necessarily disjoint), then there exists $b\in B$ such that $A\cup\{b\}\in\cal F$.
A downward closed set-system $\F$ has this property if and only if $\F$ is the family of independent sets in a matroid in which no subspace has two disjoint generators $A$ and $B$, i.e., $A\cap B=\emptyset$ and $\rk(A)=\rk(B)=\rk(A\cup B)$ is forbidden.
We do not know whether this question has been studied before.
\smallskip

\subsubsection*{Acknowledgment}

We are grateful to Bal\'azs Keszegh for his valuable remarks concerning Lemma 2.1, 
and to our anonymous reviewer for the careful reading and comments.

\end{document}